\documentclass[11pt]{article}
\usepackage{amsmath}
\usepackage{amssymb}
\usepackage{graphicx}
\usepackage{amscd}
\usepackage{amsfonts}
\usepackage{bm}                                                
\usepackage{float}
\usepackage{latexsym}
\usepackage{lscape}
\usepackage{mathrsfs}
\usepackage{verbatim}
\usepackage{color}

\newcommand{\ds}{\displaystyle}
\newtheorem{theorem}{Theorem}[section]

\newtheorem{lemma}{Lemma}[section]

\newcommand{\mcL}{\ensuremath{\mathcal{L}}}

\newcommand{\mcS}{\ensuremath{\mathcal{S}}}

\newcommand{\mrI}{\ensuremath{\mathrm{I}}}

\definecolor{armygreen}{rgb}{0.29, 0.33, 0.13}

\newcommand{\be}{\begin{equation}}
\newcommand{\ee}{\end{equation}}
\newcommand{\bea}{\begin{eqnarray}}
\newcommand{\eea}{\end{eqnarray}}
\newcommand{\beas}{\begin{eqnarray*}}
\newcommand{\eeas}{\end{eqnarray*}}
\newcommand{\ba}{\begin{array}}
\newcommand{\ea}{\end{array}}

\def\qed{\hbox{\vrule width 6pt height 6pt depth 0pt}}

\title{Numerical approximations for the variable coefficient fractional diffusion equations with non-smooth data} 
\author{
	Xiangcheng Zheng \thanks{Department of Mathematics, University of South Carolina, Columbia,
		South Carolina 29208, USA. email: {\tt xz3@math.sc.edu \& hwang@math.sc.edu}.} 
	\and
	V.J.~Ervin\thanks{Department of Mathematical Sciences,
	  Clemson University, Clemson, South Carolina 29634-0975, USA.
	  email: {\tt vjervin@clemson.edu}. }
	\and 
	 Hong Wang $^*$ } 

\date{\today}

\begin{document}

\maketitle

\begin{abstract}
In this article we study the numerical approximation of a variable coefficient
fractional diffusion equation.
 Using a change of variable, the variable coefficient fractional diffusion equation
is transformed into a constant coefficient fractional diffusion equation of the same order. The transformed
equation retains the desirable stability property of being an elliptic equation. A spectral approximation 
scheme is proposed and analyzed for the transformed equation, with error estimates for the approximated solution derived. An approximation to the unknown of the variable coefficient fractional diffusion equation is then obtained by post processing the computed approximation to the transformed equation. Error estimates are also presented for the approximation to the unknown of the variable coefficient equation with both smooth and non-smooth diffusivity coefficient and right-hand side. Three numerical experiments are given whose convergence results
are in strong agreement with the theoretically derived estimates.

\end{abstract}



\setcounter{equation}{0}
\setcounter{figure}{0}
\setcounter{table}{0}
\setcounter{theorem}{0}
\setcounter{lemma}{0}
\setcounter{corollary}{0}
\section{Introduction}
It has been shown that fractional partial differential equations (PDEs) can accurately model challenging phenomena including anomalous transport, long-range time memory and spatial interactions \cite{BenWhe,MetKla00}. Extensive research has been conducted on fractional PDEs in terms of their modeling, analysis, numerical approximations and applications. In a representative piece of work Ervin and Roop \cite{ErvRoo05} proved the wellposedness of the Galerkin weak formulation of linear elliptic space fractional diffusion equations (FDEs) of order $1 < \alpha < 2$ on the Sobolev space $H^{\alpha/2}_0$.
They also proved optimal-order error estimates of its finite element approximations in the energy and $L^{2}$ norms, assuming that the true solution and the solution to the dual problem for an $L^{2}$ right-hand side have full regularity. However, it was later realized that the smoothness of the coefficients and source term for these space fractional differential equations cannot ensure the smoothness of their solution \cite{ErvHeuRoo,JinLaz,WanYanZhu14,WanZha}. This is in sharp contrast to integer-order linear elliptic PDEs \cite{Eva,GilTru}. For this reason, the usual smoothness assumptions on the true solutions to fractional PDEs in the analysis of the numerical approximations are inappropriate.

It turns out that the spectral methods are particularly well suited for the accurate approximation of FDEs, as they provide a clean expression of the true solution to FDEs for the convenience of analysis \cite{ErvHeuRoo,MaoGeo}, and ordinarily lead to a diagonal stiffness matrices (at least for constant-coefficient FDEs). This is in contrast to the dense stiffness matrices generated from the finite element, finite difference, or finite volume approximations.
 Mao et al. \cite{MaoCheShe} analyzed the regularity of the solution to a symmetric case of the FDE and developed corresponding spectral methods. The solution structure to the general case was resolved completely in \cite{ErvHeuRoo}, in which a spectral method utilizing the weighted Jacobi polynomial was studied and a priori error estimates derived. The two-sided FDE with constant coefficient and Riemann-Liouville fractional derivative was investigated in \cite{MaoGeo}, by employing a Petrov-Galerkin projection in a properly weighted Sobolev space using two-sided Jacobi polyfractonomials as test and trial functions. In \cite{Hao17} and \cite{Hao18}, the regularity of the two-sided fractional reaction-diffusion
 and advection-diffusion-reaction equations are analyzed in the weighted Sobolev spaces, based on which the optimal (or sub-optimal) convergence rates of the spectral Galerkin or Petrov-Galerkin method are proved.

The variable diffusivity $K$ presents another bottleneck of FDEs. It was shown in \cite{WanYan} that the Galerkin weak formulation may lose its coercivity for a smooth $K(x)$ with positive lower and upper bounds, which increases the difficulties for the stability and convergence analysis and accurate simulations. To circumvent this issue, an indirect Legendre spectral Galerkin method was developed for the FDE in \cite{WanZha}, in which the high-order convergence rates of numerical approximations were proved only under the regularity assumptions of coefficients and right-hand side term. In \cite{LiCheWan}, with the introduction of an auxiliary variable, a mixed approximation was developed for an FDE and the corresponding error estimates were proved. In \cite{MaoShe}, a spectral Galerkin method for a different variable coefficient FDE was analyzed, in which the outside and inside fractional derivatives are chosen carefully so that the corresponding Galerkin weak formulation are self-adjoint and coercive. 

 Recently the wellposedness of the variable coefficient FDE (\ref{ModelV}) was investigated in \cite{ZheErvWan}, in which the existence and uniqueness of the solution to the proposed model was proven for any $f\in L^2_{\omega^{(\beta-1,\alpha-\beta-1)}}$, with the space defined by (\ref{norm}). A spectral approximation method was then studied and several error estimates were derived based on the regularity of the right-hand side term.  
In this paper we continue to investigate model (\ref{ModelV}) using a different approach than that used in \cite{ZheErvWan}. We prove in this paper that the model is wellposed for $f$ belonging to a larger space $ L^2_{\omega^{(\beta,\alpha-\beta)}}$, which extends the wellposedness results in \cite{ZheErvWan}. A spectral approximation scheme is proposed and the error estimates are proved to be dependent on the weaker norms of $f$ without loss of accuracy. In addition, we follow the idea of the K-method of interpolation to determine the range of the index of the weighted Sobolev space that the power function belongs to, which provides the theoretical support for estimating the convergence rate of the proposed method. 

This paper is organized as follows. In Section 2 we present the formulation of the model and introduce notation and key lemmas used in the analysis. The wellposedness of the model and the regularity of its solution are studied in Section 3, based on which the spectral approximation method is formulated and a detailed analysis of its convergence is proved. Three numerical experiments are presented in Section 4 whose results demonstrate the sharpness of the derived error estimates.

\setcounter{equation}{0}
\setcounter{figure}{0}
\setcounter{table}{0}
\setcounter{theorem}{0}
\setcounter{lemma}{0}
\setcounter{corollary}{0}
\section{Model problem and preliminaries}
In this paper we consider the following homogeneous Dirichlet boundary-value problem of a two-sided Caputo flux FDE, 
which is obtained by incorporating a fractional Fick's law into a conventional local mass balance law 
 with a variable fractional diffusivity \cite{CasCar,ZhaBen}:
\begin{align}
 -D \big (  \big (r \,{}_0I_x^{2-\alpha} +  (1-r) \,  {}_xI_1^{2-\alpha} \big ) \, K(x) \,  Du(x) \big ) &= f(x), 
 ~~x \in  \mrI \, := (0,1),    \label{ModelV}  \\
u(0) = u(1) = 0.  \label{BC1}
\end{align}
Here $1 < \alpha < 2$, $D$ is the first-order differential operator, $K(x)$ is the fractional diffusivity with 
$0 < K_{m} \leq K(x) \leq K_{M} < \infty$, $0\leq r\leq 1$ indicates the relative weight of forward versus backward 
transition probability and $f(x)$ the source or sink term. 
The left and right fractional integrals of order $0 < \sigma < 1$ are defined as \cite{Pod}
\begin{equation*}
{}_0I_x^{\sigma}w(x) := \frac{1}{\Gamma(\sigma)} \int_0^x w(s)(x-s)^{\sigma-1}ds, \quad {}_xI^{\sigma}_1w(x) := \frac{1}{\Gamma(\sigma)} \int_x^1 w(s)(s-x)^{\sigma-1} \, ds \, ,
\end{equation*}
where $\Gamma(\cdot)$ is the Gamma function. 

We introduce notation and properties used subsequently in our discussion of the approximation scheme and
in its error analysis. 

Let $\Omega$ be a bounded open interval and $\omega(x) > 0, ~x \in \Omega$ be a smooth function. We define the weighted $L^2$ space,  $L_{\omega}^{2}(\Omega)$, and $L^{2}$ weighted inner product as
\begin{equation}\label{norm} L_{\omega}^{2}(\Omega) \, := \, \bigg \{ f(x) \, : \, \| f\|_{\omega}^2 := 
\int_{\Omega} \omega(x) \, f(x)^{2} \, dx \ < \ \infty \bigg \}, \quad
(f \, , \, g)_{\omega} \, := \, \int_{\Omega} \omega(x) \, f(x) \, g(x) \, dx \, . \end{equation}

In addition, let
$\mathbb{N}_{0}  := \mathbb{N} \cup {0}$,
and  $\omega^{(\alpha,\beta)}$ be a weighting function defined on $\Omega$ and indexed by $\alpha$ and $\beta$. 
For any $m \in \mathbb{N}$, we introduce the following weighted Sobolev spaces \cite{GuoWan,SheTan}
$$\ds H^m_{\omega^{(\alpha,\beta)} }(\Omega) := \bigg \{ v \, : \| v \|_{m,\omega^{(\alpha,\beta)}}^2 := \sum_{j=0}^m \big | v \big |_{j,\omega^{(\alpha,\beta)}}^2 = \sum_{j=0}^m \big \| D^j v \big \|_{\omega^{(\alpha+j,\beta+j)}}^2 < \infty \bigg \} .$$  
For $t \in \mathbb{R}^{+} \backslash \mathbb{N}_{0}$, $H^{t}_{\omega^{(a \, , \, b)}}(0,1)$ 
is defined by the $K$-method of interpolation, and
for $t \in \mathbb{R}^{-}$, $H^{t}_{\omega^{(a \, , \, b)}}(0,1)$ is defined by duality.

The Jacobi polynomials $P_{n}^{(\alpha , \beta)}(x)$ are defined by \cite{SheTan,Sze}
\begin{equation}\label{Pn}
P_{n}^{(\alpha , \beta)}(x) := \sum_{k = 0}^{n} p_{n,k} (x - 1)^{n - k} (x + 1)^k, ~~x \in (-1,1), \quad p_{n,k} :=  \frac{1}{2^{n}} {n + \alpha \choose{k}}{n + \beta \choose{n-k}}.
\end{equation}
Let $G_{n}^{(\alpha,\beta)}(x)$ denote the translated and dilated Jacobi polynomials to the interval $[0,1]$:
\begin{equation}\label{Gn}
G_{n}^{(\alpha,\beta)}(x) := P_{n}^{(\alpha,\beta)}(2x-1), ~~x \in [0,1],
\end{equation}
and summarize the properties of $G_{n}^{\alpha,\beta}$ in the following lemma

\begin{lemma}\label{lem:Gn}
For $\omega^{(\alpha , \beta)}(x) = (1 - x)^{\alpha} x^{\beta}$,
the polynomials $G_n^{(\alpha,\beta)}$ have the following orthogonality and norm properties 
\begin{equation}\label{GnPerp}
\int_0^1  \omega^{(\alpha,\beta)}(x) G_{j}^{(\alpha , \beta)}(x) \, G_{k}^{(\alpha , \beta)}(x)  dx = \delta_{j,k} \,  |\| G_{j}^{(\alpha , \beta)}  |\|^{2}, \qquad \omega^{(\alpha , \beta)}(x) := (1 - x)^{\alpha} x^{\beta},
\end{equation}
where $\delta_{j,k} = 1$ if $j \neq k$ and 0 otherwise, and
\begin{equation}\label{Gnorm}
 |\| G_{j}^{(\alpha , \beta)}  |\| := \Big ( \frac{1}{(2j \, + \, \alpha \, + \, \beta \, + 1)}    \frac{\Gamma(j + \alpha + 1) \, \Gamma(j + \beta + 1)}{\Gamma(j + 1) \, \Gamma(j + \alpha + \beta + 1)}\Bigr)^{1/2} =  |\| G_{j}^{(\beta,\alpha)}  |\| .
\end{equation}
In addition, $G_n^{(\alpha,\beta)}$ satisfies
\begin{equation}\label{DGn}\begin{array}{l}
\ds D^{k}G_{n}^{(\alpha,\beta)}(x) = \frac{\Gamma(n + k + \alpha + \beta + 1)}{\Gamma(n + \alpha + \beta + 1)} 
G_{n - k}^{(\alpha + k,\beta + k)}(x), \quad 0 \le k \le n;\\[0.15in]
\ds D^k \Big ( \omega^{(\alpha + k,\beta + k)}(x) G_{n-k}^{(\alpha + k,\beta + k)}(x) \Big ) =  
\frac{(-1)^{k} \, n!}{(n - k)!} \, \omega^{(\alpha,\beta)}(x)\,G_{n}^{(\alpha,\beta)}(x), \quad 0 \le k \le n. 
\end{array}\end{equation}   
Finally, $G_n^{(\alpha,\beta)}$ have the following norm relation
\begin{equation}\label{GnNormRelation}
\frac{1}{2} \, \le \, \frac{ |\| G_{j}^{(\alpha - \beta \, , \, \beta)} |\|^{2}}{ |\| G_{j + 1}^{(\beta - 1 \, , \, \alpha - \beta - 1)} |\|^{2}} 
	\, = \, \frac{j + 1}{j + \alpha} \, \le \, 1, \, \quad j \ge 0.
	\end{equation}	
\end{lemma}

{\it Proof.} The orthogonality property (\ref{GnPerp}) of $G_{n}^{(\alpha,\beta)}$ is a direct consequence of the orthogonality relation of $P_{n}^{(\alpha,\beta)}$ \cite{SheTan,Sze}
$$\begin{array}{c}
\ds \int_{-1}^1 \tilde{\omega}^{(\alpha,\beta)}(x) P_j^{(\alpha , \beta)}(x) P_k^{(\alpha,\beta)}(x) dx = \delta_{j,k} \big \| P_{j}^{(\alpha,\beta)} \big \|^2_{\tilde{\omega}^{(\alpha,\beta)}}, 
\ \ \mbox{ where } \ \  \tilde{\omega}^{(\alpha,\beta)}(x) := (1 - x)^{\alpha} (1 + x)^{\beta}, \\[0.1in]
 \mbox{ Using } \ \ \big \| P_{j}^{(\alpha,\beta)} \big \|_{\tilde{\omega}^{(\alpha,\beta)}} = \Big( \frac{2^{(\alpha + \beta + 1)}}{(2j \, + \, \alpha \, + \, \beta \, + 1)} \frac{\Gamma(j + \alpha + 1) \, \Gamma(j + \beta + 1)}{\Gamma(j + 1) \, \Gamma(j + \alpha + \beta + 1)}\Big)^{1/2} =  \big \| P_{j}^{(\beta,\alpha)} \big \|_{\tilde{\omega}^{(\beta,\alpha)}}
\end{array}$$
and the following relation between $G_{n}^{(\alpha,\beta)}$ and $P_{n}^{(\alpha,\beta)}$ leads to \eqref{Gnorm}.
$$ \int_{-1}^{1} \tilde{\omega}^{(\alpha,\beta)}(x) \, P_{j}^{(\alpha , \beta)}(x) P_{k}^{(\alpha , \beta)}(x) dx
= 2^{\alpha + \beta + 1} \int_0^1  \omega^{(\alpha,\beta)}(x) G_{j}^{(\alpha , \beta)}(x) \, G_{k}^{(\alpha , \beta)}(x)dx.$$

The two equations in (\ref{DGn}) are direct consequences of the following equations for $P_{n}^{(\alpha,\beta)}(x)$ \cite[equations (2.15) and (2.19)]{MaoCheShe}  
\begin{equation*}\begin{array}{l}
\ds D^k P_{n}^{(\alpha,\beta)}(x) = \frac{\Gamma(n + k + \alpha + \beta + 1)}{2^{k} \, \Gamma(n + \alpha + \beta + 1)} P_{n - k}^{(\alpha + k,\beta + k)}(x), \\[0.15in]
\ds D^k \Big (\tilde{\omega}^{\alpha + k,\beta + k}(x) P_{n - k}^{(\alpha + k,\beta + k)}(x) \Big )
= \frac{(-1)^{k} \, 2^{k} \, n!}{(n - k)!} \, \tilde{\omega}^{\alpha,\beta}(x) P_{n}^{(\alpha,\beta)}(x).
\end{array}\end{equation*}
The norm relation (\ref{GnNormRelation}) is derived in \cite{ZheErvWan}.
\qed

Let $\mcS_{N}$ denote the space of polynomials of degree $\le N$. We define the weighted $L^2$ orthogonal projection
$P_{N, a, b} : \, L^{2}_{\omega^{(a , b)}}(0,1) \rightarrow \mcS_{N}$ by the condition
\begin{equation}\label{Proj}
\big ( v \, - \,  P_{N, a, b}v \ , \ \phi_N \big)_{\omega^{(a , b)}} \ = \ 0 \, , \ \ \forall \phi_N \in \mcS_{N}.
\end{equation}
\begin{lemma}\label{lem:Approx} \cite[Theorem 2.1]{GuoWan}
For $\mu \in \mathbb{N}_{0}$ and $v \in H^{t}_{\omega^{(a \, , \, b)}}(0,1)$, with $0 \le \mu \le t$, there exists a
constant $C$, independent of $N, \, \alpha$ and $\beta$ such that
\begin{equation}\label{Approx}
\big \| v  -  P_{N, a, b} v \|_{\mu , \omega^{(a , b)}} \ \le \ C \, 
\big ( N \, (N + a + b) \big )^{\frac{\mu - t}{2}} \, | v |_{t ,\omega^{(a , b)}}.
\end{equation}
\end{lemma}  

\textbf{Remark}: In \cite{GuoWan} \eqref{Approx} is stated for $t \in \mathbb{N}_{0}$. The result extends to 
$t \in \mathbb{R}^{+}$ using interpolation.

\setcounter{equation}{0}
\setcounter{figure}{0}
\setcounter{table}{0}
\setcounter{theorem}{0}
\setcounter{lemma}{0}
\setcounter{corollary}{0}
\section{Approximation scheme}
\label{secapx}

\subsection{Motivation for the approximation scheme}

Introduce $D^{-1} \, : \, L^{1}(\mrI) \longrightarrow H^{1}(\mrI)$, defined by $D^{-1} g(x) \, := \, \int_{0}^{x} \, g(s) \, ds$.

Rewrite \eqref{ModelV} as 
\[
-D \big (  \big (r \,{}_0I_x^{2-\alpha} +  (1-r) \,  {}_xI_1^{2-\alpha} \big ) \, D 
\underbrace{ D^{-1} \, K(x) \,  Du(x)}_{:= \, \widetilde{w}(x)} \big )  \ =  \ f(x) \, .
\]

Consider $\widetilde{w}(x) \ = \ D^{-1} \, K(x) \,  Du(x) \ = \ \int_{0}^{x} K(s) \, D \, u(s) \, ds$. Note that
$\widetilde{w}(0) = 0$. Then,
\begin{align}
D \, \widetilde{w}(x) &= \ K(x) \, D \, u(x)   \nonumber  \\
\Longrightarrow \ \ \ u(x) &= \ \int_{0}^{x} \frac{D \, \widetilde{w}(s)}{K(s)} \, ds \, .  \label{htfu1} \\
\mbox{Now, } \  \  \ u(1) \, = \, 0 \  \ \Longrightarrow \ \ 
 \int_{0}^{1} \frac{D \, \widetilde{w}(s)}{K(s)} \, ds &= \ 0 \, .  \label{bcd1}
\end{align} 

Hence if we could determine $\widetilde{w}(x)$ such that
\begin{align*}
\mcL_{r}^{\alpha} \widetilde{w}(x) \ := \  -D \big (  \big (r \,{}_0I_x^{2-\alpha} +  (1-r) \,  {}_xI_1^{2-\alpha} \big ) \, 
D\widetilde{w}(x) \big ) &= f(x),  ~~x \in  \mrI ,      \\
\mbox{subject to } \ \widetilde{w}(0) \, = \, 0 \ \ \mbox{ and } \int_{0}^{x} \frac{D \widetilde{w}(s)}{K(s)} \, ds &= \, 0 \, ,
\end{align*}
then our solution to \eqref{ModelV},\eqref{BC1} would be given by \eqref{htfu1}.  With this in mind, 
consider the problem: 
\textit{Determine $w(x)$ satisfying}
\begin{align}
\mcL_{r}^{\alpha} w(x)  &= f(x), 
 ~~x \in  \mrI ,    \label{ModelVI}  \\
\mbox{subject to } \ \ w(0) &= w(1) = 0.  \label{BC2}
\end{align}
Let $\beta\in [\alpha-1, 1]$ be determined by $(1-r)\sin(\pi \beta)=r\sin(\pi(\alpha-\beta))$. The following theorem ensures the well-posedness of this problem.
\begin{theorem}[\cite{ErvHeuRoo, JiaChenErvin}]   \label{thmWellposedw}
	Let $f (x) \in L^{2}_{\omega^{(\beta , \alpha - \beta)}}(\mrI)$. Then, there exists a unique solution
	$w(x) 
	\in L^{2}_{\omega^{(-(\alpha - \beta) , -\beta)}}(\mrI)$  satisfying \eqref{ModelVI},\eqref{BC2}.
	In addition, there exists $C > 0$ such that
	\begin{equation}
	\| w \|_{\omega^{(-(\alpha - \beta) , -\beta)}} \, + \, \| D w \|_{\omega^{(-(\alpha - \beta) + 1 \,  , \, -\beta + 1)}} 
	\ \le \ C \, \| f \|_{\omega^{(\beta , \alpha - \beta)}}
	\end{equation}
\end{theorem}
\mbox{ } \hfill \qed

From \cite{ErvHeuRoo} we have that  $k_{1}(x) \, := \, \int_{0}^{x} (1 - s)^{\alpha - \beta - 1} \, s^{\beta - 1} \, ds \, \in \, 
ker(\mcL_{r}^{\alpha}).$

\begin{equation}
\mbox{Let } \ \ \widetilde{w}(x) \ = \ C_{1} \, k_{1}(x) \ + \ w(x) . 
\label{defwt}
\end{equation}

We have that $\mcL_{r}^{\alpha} \widetilde{w}(x) \ = \ f(x)$, and $\widetilde{w}(0) \, = \, 0$.

The condition \eqref{bcd1} combined with \eqref{defwt} implies
\begin{align} 
 \int_{0}^{1} \frac{1}{K(s)} \big( C_{1} \, (1 - s)^{\alpha - \beta -1} \, s^{\beta - 1}  &+ \ D w(s) \big) \, ds \ = \ 0 \, ,  \nonumber \\
 \Longrightarrow \ \ C_{1} &= \frac{- \int_{0}^{1} \frac{D w(s)}{K(s)} \, ds}%
{\int_{0}^{1} \frac{(1 - s)^{\alpha - \beta -1} \, s^{\beta - 1}}{K(s)} \, ds}  \label{f4c1}  \\
\mbox{ (if 
$D(\frac{1}{K}) \in L^{2}_{\omega^{(\alpha - \beta \, , \, \beta)}}$
)} \ \ \ 
 \quad \quad \quad  \quad 
&= \frac{\int_{0}^{1} w(s) \, \frac{d}{ds} \left( \frac{1}{K(s)} \right) \, ds}%
{\int_{0}^{1} \frac{(1 - s)^{\alpha - \beta -1} \, s^{\beta - 1}}{K(s)} \, ds}  \nonumber  \\
 \quad \quad \quad  \quad 
&= \frac{\int_{0}^{1} w(s) \,  \frac{K'(s)}{K(s)^{2}}  \, ds}%
{\int_{0}^{1} \frac{(1 - s)^{\alpha - \beta -1} \, s^{\beta - 1}}{K(s)} \, ds}  \, . \label{f4c2}  
\end{align}

Let 
\begin{align}
\ den \, := \, \int_{0}^{1}  \frac{(1 - s)^{\alpha - \beta -1} \, s^{\beta - 1}}{K(s)} \, ds \, , \quad
c_{1} \, := \, \int_{0}^{1} \frac{- D w(s)}{K(s)} ds \, . \label{eq19}
\end{align}

Then $C_1=c_1/den$ can be bounded by
\begin{align}
|C_1| &= \frac{1}{den}\ 
\big| \int_{0}^{1} D  w (s) \, \frac{1}{K(s)} \, ds  \big| 
\ \le \frac{1}{den}\int_{0}^{1} \big| D  w (s) \big| \, \frac{1}{K(s)} \, ds   \label{bndC0}  \\
&=  \frac{1}{den}\ 
\int_{0}^{1} \omega^{( (-(\alpha - \beta) + 1)/2 \, , \, (-\beta + 1)/2)}(s) \, \big| D  w (s) \big| \, 
\omega^{( (\alpha - \beta - 1)/2 \, , \, (\beta - 1)/2)}(s) \, \frac{1}{K(s)} \, ds  \nonumber \\
&\leq \frac{1}{den}\ 
\| D  w  \|_{\omega^{( -(\alpha - \beta) + 1 \, , \, -\beta + 1)}} \, 
\bigg\| \frac{1}{K(s)} \bigg\|_{\omega^{( (\alpha - \beta - 1) \, , \, (\beta - 1))}}    \nonumber   \\
&\le C   \, \| f \|_{\omega^{(\beta  \, , \, \alpha - \beta )}}  
\ \ \mbox{(using \mbox{Theorem} ~\ref{thmWellposedw})} \, .  \label{bndC1}
\end{align}

Combining \eqref{htfu1},\eqref{defwt}, \eqref{f4c1}, (\ref{bndC1}) and Theorem \ref{thmWellposedw} we have the following.
\begin{theorem} \label{exun1}
For $f(x) \in L^{2}_{\omega^{(\beta , \alpha - \beta)}}(\mrI)$ and there exists a
unique solution $u(x) \in L^{\infty}(\mrI)$ to \eqref{ModelV},\eqref{BC1}, given by
\begin{equation}
u(x) \ = \ C_{1} \int_{0}^{x}  \frac{(1 - s)^{\alpha - \beta -1} \, s^{\beta - 1}}{K(s)} \, ds \ + \ 
   \int_{0}^{x}  \frac{D w(s)}{K(s)} \, ds \, ,
\label{ans4u}
\end{equation}
where $w(x)$ is determined by \eqref{ModelVI},\eqref{BC2} and $C_{1}$ by \eqref{f4c1}.

Additionally, for $\epsilon_{1}, \epsilon_{2} > 0$ there exists $C > 0$ such that
\begin{equation}
\| u \|_{L^{\infty}} \ + \ \| u \|_{\omega^{(-1 + \epsilon_{1} \, , \, -1 + \epsilon_{2})}} \ \le \ C \,  \| f \|_{\omega^{(\beta  \, , \, \alpha - \beta )}} \, .
\label{uest11}
\end{equation}

\end{theorem}
\textbf{Proof}: \\
It is straightforward to show that $u(x)$ given by \eqref{ans4u} satisfies \eqref{ModelV},\eqref{BC1}. Next we
show that there exists a unique solution to \eqref{ModelV},\eqref{BC1}.

Assume that $u_{1}(x)$ and $u_{2}(x)$ are solutions of \eqref{ModelV},\eqref{BC1}. Let 
$z_{1}(x)$ and $z_{2}(x)$ be defined by:
\begin{align}
    z_{1}(x) \ = \ \int_{0}^{x} K(s) \, D u_{1}(s) \, ds \, , & \quad z_{1}(x) \ = \ \int_{0}^{x} K(s) \, D u_{1}(s) \, ds \, \, ,
    \nonumber \\
\mbox{i.e., } \ \ D z_{1}(x) \ = \ K(x) \, D u_{1}(x)   \, ,  \ \ \  z_{1}(0) \, = \, 0 \, , & \mbox{  and  }  \ \ 
    D z_{2}(x) \ = \ K(x) \, D u_{2}(x)   \, ,  \ \ \  z_{2}(0) \, = \, 0 \, .  \label{defz12}
\end{align}

Note that $\mcL_{r}^{\alpha} ( z_{1} \, - \, z_{2}) \ = \ 0$. Hence, $( z_{1} \, - \, z_{2}) \in ker(\mcL_{r}^{\alpha} )$.
Thus from \cite{ErvHeuRoo}, for constants $A$ and $B$,
\begin{align}
( z_{1} \, - \, z_{2})(x) \ = \ A &+ \ B \, \int_{0}^{x} (1 - s)^{\alpha - \beta - 1} \, s^{\beta - 1} \, ds \, .  \nonumber \\
\mbox{As } \ ( z_{1} \, - \, z_{2})(0) \, = \, 0 \, \ \Longrightarrow \ A &= \ 0 \, .  \nonumber  \\
\mbox{Then,  } \quad  D ( z_{1} \, - \, z_{2})(x) &= \ B \,  (1 - x)^{\alpha - \beta - 1} \, x^{\beta - 1}  \nonumber \\
  \Longrightarrow \ \ D ( u_{1} \, - \, u_{2})(x) &= \ B \,  \frac{1}{K(x)} \, (1 - x)^{\alpha - \beta - 1} \, x^{\beta - 1} \nonumber \\
  \Longrightarrow \ \  \mbox{ (using $(u_{1} \, - \, u_{2})(0) \, = \, 0$) } \ \ \ 
  ( u_{1} \, - \, u_{2})(x) &= \ B \,  \int_{0}^{x} \, \frac{1}{K(s)} \, (1 - s)^{\alpha - \beta - 1} \, s^{\beta - 1} \, ds  \, . \nonumber \\
  \mbox{ As the integrand is nonnegative, $(u_{1} \, - \, u_{2})(1) \, = \, 0$ }  \Longrightarrow \ \  
  B &= 0 \, , \nonumber \\
  \Longrightarrow \ \  u_{1}(x) &= \ u_{2}(x) \, .  \nonumber 
\end{align}

Using \eqref{bndC0} and (\ref{bndC1}) 
\begin{align}
| u(x)  | &= \ \big| C_{1}  \int_{0}^{x}  \frac{(1 - s)^{\alpha - \beta -1} \, s^{\beta - 1}}{K(s)} \, ds 
\ + \   \int_{0}^{x}  \frac{D w (s)}{K(s)} \, ds  \, \big| \, ,   \nonumber \\
&\le \ C \, \big| C_{1} \big| \ + \  \int_{0}^{1} \big| D w (s) \big| \, \frac{1}{K(s)} \, ds  \leq C  \, \| f \|_{\omega^{(\beta  \, , \, \alpha - \beta )}}  \, .  \label{uinfb} 
\end{align}

Consequently, we obtain
\begin{align}
\|u\|^2_{\omega^{(-1 + \epsilon_{1} \,  , \, -1 + \epsilon_{2})}}&=\int_0^1 \, u^2(x)\, 
\omega^{(-1 + \epsilon_{1} \,  , \, -1 + \epsilon_{2})}(x) \, dx  \nonumber \\
&\leq C  \, \| f \|^2_{\omega^{(\beta  \, , \, \alpha - \beta )}}\int_0^1 \, \omega^{(-1 + \epsilon_{1} \,  , \, -1 + \epsilon_{2})}(x) \, dx\leq C  \, \| f \|^2_{\omega^{(\beta  \, , \, \alpha - \beta )}} \, .  \label{uombd}
\end{align}
Estimate \eqref{uest11} then follows from \eqref{uinfb} and \eqref{uombd}. \\
\mbox{  } \hfill \qed

\subsection{Approximation scheme}
\label{ssec_bvp2}
To compute an approximation to $u(x)$, $u_{N}(x)$, we firstly compute an approximation to
$w(x)$, $w_{N}(x)$, satisfying \eqref{ModelVI},\eqref{BC2}, and then use $w_{N}(x)$ in place
of $w(x)$ in \eqref{f4c1} and \eqref{ans4u} to obtain $u_{N}(x)$.

\subsubsection{Approximation of $w(x)$ satisfying \eqref{ModelVI},\eqref{BC2}}
\label{sssec_w}
Proceeding as in \cite{ErvHeuRoo}, 
 $f(x) \in L^{2}_{\omega^{(\beta , \alpha - \beta)}}(\mrI)$ may be expressed as 
$f(x) \ = \ \sum_{i = 0}^{\infty} \frac{f_{i}}{ |\| G_{i}^{(\beta , \alpha - \beta)} |\|^{2}} \, G_{i}^{(\beta , \alpha - \beta)}(x)$, 
where $f_{i}$ is given by
\begin{equation}
f_{i} \, := \, \int_{0}^{1} \,  \omega^{(\beta , \alpha - \beta)}(x) \, f(x) \, G_{i}^{(\beta , \alpha - \beta)}(x)  \, dx \, .
\label{deffis}
\end{equation}

With $f_{i}$ defined in \eqref{deffis}, let 
\begin{equation}
f_{N}(x) \ = \ \sum_{i = 0}^{N} \frac{f_{i}}{ |\| G_{i}^{(\beta , \alpha - \beta)} |\|^{2}} \, G_{i}^{(\beta , \alpha - \beta)}(x) 
 \  \mbox{ and } 
w_{N}(x) \ = \ \omega^{( \alpha - \beta , \beta )}(x) \, \sum_{i = 0}^{N} c_{i} \, G_{i}^{( \alpha - \beta , \beta )}(x) \, , 
 \label{defuNs} 
\end{equation}
\begin{equation}
\mbox{ where} \ 
\lambda_{i} \ = \  \frac{- \, \sin(\pi \alpha)}{\sin(\pi (\alpha - \beta)) \, + \, \sin(\pi \beta)} \,
\frac{\Gamma(i + 1 + \alpha)}{\Gamma(i + 1)}  \  \mbox{ and }
c_{i}  \ = \ \frac{1}{ \lambda_{i} \, |\| G_{i}^{( \beta , \alpha - \beta  )} |\|^{2}  }f_{i}  \, .
\label{defuNs2}
\end{equation}

Using Stirling's formula we have that 
\begin{equation}
\lim_{n \rightarrow \infty} \frac{\Gamma(n + \mu)}{\Gamma(n) \, n^{\mu}} \ = \ 1 \, ,  \mbox{ for } \mu  \in \mathbb{R}.
\ \mbox{Thus } \lambda_{i} \sim (i + 1)^{\alpha} \, .
\label{lbalge}
\end{equation}

\begin{theorem}[\cite{ErvHeuRoo, JiaChenErvin}]   \label{thmuspg}
Let $f (x) \in L^{2}_{\omega^{(\beta , \alpha - \beta)}}(\mrI)$ and $w_{N}(x)$ be as defined in \eqref{defuNs}. Then,
$w(x) \ := \ \lim_{N \rightarrow \infty} w_{N}(x) \ = \ \omega^{( \alpha - \beta , \beta )}(x) \, \sum_{j = 0}^{\infty}  c_{j} \,
G_{j}^{( \alpha - \beta , \beta )}(x) \, 
\in L^{2}_{\omega^{(-(\alpha - \beta) , -\beta)}}(\mrI)$ and satisfies \eqref{ModelVI},\eqref{BC2}.
\end{theorem}
\mbox{ } \hfill \qed

\begin{theorem} 	\label{thmAP2}
	For $f(x) \in H^{t}_{\omega^{(\beta  \, , \, \alpha - \beta)}}(\mrI)$,
$t \ge 0$,	
	and $w_{N}(x)$ given by  \eqref{defuNs}, there exists $C > 0$ such that
	\begin{align}
	\| w \, - \, w_{N} \|_{\omega^{(-(\alpha - \beta) \, , - \beta)}} \ &\le\, C \, (N + 2)^{- \alpha } (N(N+\alpha))^{-t \,/ \,2}   \, \| f \|_{t,\,\omega^{(\beta  \, , \, \alpha - \beta )}} , \ \mbox{ and}
	\label{apertg2}\\
		\|D( w \, - \, w_{N}) \|_{\omega^{(-(\alpha - \beta)+1 \, , -\beta + 1))}} \ &\le\, C \, (N + 2)^{- (\alpha-1) } (N(N+\alpha))^{-t \,/ \,2}   \, \| f \|_{t,\,\omega^{(\beta  \, , \, \alpha - \beta )}} . \label{error:Dw}
	\end{align} 
\end{theorem} 
\textbf{Proof}: Using the definition of the $ \| \cdot \|_{\omega^{(-(\alpha - \beta) \, , - \beta)}} $ norm, 
\begin{align*}
\| w \, - \, w_{N} \|_{\omega^{(-(\alpha - \beta) \, , - \beta)}}^{2} &=
\int_{0}^{1} \omega^{(-(\alpha - \beta) \, , - \beta)}(x) \bigg( \omega^{(\alpha - \beta \, , \, \beta)}(x) \, \sum_{i = N+1}^{\infty}  
\frac{G^{(\alpha - \beta \, , \, \beta)}_{i}(x)}{ ( \lambda_{i} \,  |\| G^{(\beta  \, , \, \alpha - \beta )}_{i} |\|^{2} )}
\,   f_{i}  \bigg)^{2} \, dx  \\
&\le \ \max_{i \, \ge \, N+1} \left( \frac{1}{  \lambda_{i}^{2}  } \right)  \ 
\sum_{i \, = \, N +1}^{\infty}   \frac{ f^{2}_{i} }{|\| G^{(\beta \, , \, \alpha - \beta )}_{i} |\|^{2} }
\\
&= \   \frac{1}{  \lambda_{N+1}^{2}  }   \ 
\int_{0}^{1} \, \omega^{(\beta  \, , \, \alpha - \beta )}(x) \,
\left( \sum_{i \, = \, N+1}^{\infty}  \frac{ G^{(\beta  \, , \, \alpha - \beta )}_{i}(x) }{|\| G^{(\beta  \, , \, \alpha - \beta )}_{i} |\|^{2} } \, 
f_{i}
\right)^{2} \, dx  \\
&= \ \frac{1}{  \lambda_{N+1}^{2}  }   \
\int_{0}^{1} \, \omega^{(\beta  \, , \, \alpha - \beta )}(x) \, (f(x)-P_{N,\beta,\alpha-\beta}f(x))^{2} \, dx    \\[0.05in]
&\le \ C \, (N + 2)^{-2 \, \alpha}   \ 
\| f-P_{N,\beta,\alpha-\beta}f \|_{\,\omega^{(\beta  \, , \, \alpha - \beta )}}^{2} \, , \quad 
( \mbox{using} ~ |\lambda_{N+1}| \sim (N+2)^\alpha) \\[0.05in]
&\ds \leq \ C \, (N + 2)^{-2 \, \alpha}   (N(N+\alpha))^{-t} 
\|f\|_{t,\,\omega^{(\beta  \, , \, \alpha - \beta )}}^{2} , \mbox{ (using \eqref{Approx})}.
\end{align*}
Similarly, using (\ref{DGn}) 
\begin{align}
&\| D(w \, - \, w_{N}) \|_{\omega^{(-(\alpha - \beta)+1 \, , -\beta + 1)}}^{2}    \nonumber \\
&~~~~=
\int_{0}^{1} \omega^{(-(\alpha - \beta)+1 \, , -\beta + 1)}(x) \bigg(D\bigg( \omega^{(\alpha - \beta \, , \, \beta)}(x) \, \sum_{i = N+1}^{\infty}  
\frac{G^{(\alpha - \beta \, , \, \beta)}_{i}(x)}{ ( \lambda_{i} \,  |\| G^{(\beta  \, , \, \alpha - \beta )}_{i} |\|^{2} )}
\,   f_{i}  \bigg)\bigg)^{2} \, dx  \nonumber \\
&~~~~=
\int_{0}^{1} \omega^{(-(\alpha - \beta)+1 \, , -\beta + 1)}(x) \bigg( \omega^{(\alpha - \beta-1 \, , \, \beta-1)}(x) \, \sum_{i = N+1}^{\infty}  
\frac{(i + 1) \, G^{(\alpha - \beta -1\, , \, \beta-1)}_{i+1}(x)}{ ( \lambda_{i} \,  |\| G^{(\beta  \, , \, \alpha - \beta )}_{i} |\|^{2} )}
\,   f_{i}  \bigg)^{2} \, dx    \nonumber \\
&~~~~= 
\sum_{i \, = \, N +1}^{\infty} \frac{(i + 1)^2}{\lambda_i^2}  \frac{ f^{2}_{i} }{|\| G^{(\beta \, , \, \alpha - \beta )}_{i} |\|^{4} }
|\| G^{(\alpha - \beta-1 \, , \, \beta-1 )}_{i+1} |\|^{2}   \nonumber \\
&~~~~= 
\sum_{i \, = \, N +1}^{\infty} \frac{(i + 1)^2}{\lambda_i^2} 
\left( \frac{i + \alpha}{i + 1} \right)^{2} \frac{ f^{2}_{i} }{|\| G^{(\beta \, , \, \alpha - \beta )}_{i} |\|^{2} }   \ \
\mbox{(using \eqref{GnNormRelation} and \eqref{Gnorm})} \, .  \label{Desthw1}
\end{align}
Using \eqref{lbalge},
\begin{equation}
 \frac{(i + \alpha)^2}{\lambda_i^2} \ \sim \ (i + \alpha)^{2} \, (i + 1)^{-2 \alpha} \ \sim \ (i + 1)^{-2 (\alpha - 1)} \, .
 \label{lbest1}
\end{equation}
Combining  \eqref{Desthw1} and \eqref{lbest1},
\begin{align*}
\| D(w \, - \, w_{N}) \|_{\omega^{(-(\alpha - \beta)+1 \, , -\beta + 1)}}^{2} 
&\le \   \frac{C}{  (N+2)^{2(\alpha-1)}  }   \ 
\sum_{i \, = \, N +1}^{\infty}  \frac{ f^{2}_{i} }{|\| G^{(\beta \, , \, \alpha - \beta )}_{i} |\|^{2} }   \\
&= \   \frac{C}{  (N+2)^{2(\alpha-1)}  }   \ 
\int_{0}^{1} \, \omega^{(\beta  \, , \, \alpha - \beta )}(x) \,
\left( \sum_{i \, = \, N+1}^{\infty}  \frac{ G^{(\beta  \, , \, \alpha - \beta )}_{i}(x) }{|\| G^{(\beta  \, , \, \alpha - \beta )}_{i} |\|^{2} } \, 
f_{i}
\right)^{2} \, dx  \\
&= \  \frac{C}{  (N+2)^{2(\alpha-1)}  }   \
\int_{0}^{1} \, \omega^{(\beta  \, , \, \alpha - \beta )}(x) \, (f(x)-P_{N,\beta,\alpha-\beta}f(x))^{2} \, dx    \\[0.05in]
&= \ C \, (N + 2)^{-2 \,( \alpha-1)}   \ 
\| f-P_{N,\beta,\alpha-\beta}f \|_{\,\omega^{(\beta  \, , \, \alpha - \beta )}}^{2} \, , \ \\[0.05in]
&\ds \leq \ C \, (N + 2)^{-2 \, (\alpha-1)}   (N(N+\alpha))^{-t} 
\|f\|_{t,\,\omega^{(\beta  \, , \, \alpha - \beta )}}^{2}.
\end{align*}
\mbox{ } \hfill \qed

\subsubsection{Approximation of $u(x)$ satisfying \eqref{ans4u}}
\label{sssec_u}
The approximation $u_{N}(x)$ of $u(x)$ is obtained by substituting $w_{N}(x)$ in place of $w(x)$ in 
\eqref{f4c1} and \eqref{ans4u}. With $den$ defined in \eqref{eq19}, let
\[
c_{1 , N} \, := \, \int_{0}^{1} \frac{- D w_{N}(s)}{K(s)} ds~~\mbox{and}~~ C_{1 , N} \, := \,  c_{1 , N} / den .
\]

Note that $| C_{1} \, - \, C_{1 , N} | \ = \ | c_{1} \, - \, c_{1 , N} | / den$. Hence the rate of convergence as $N \rightarrow \infty$
of $| C_{1} \, - \, C_{1 , N} |$ is equal to the rate of convergence of $| c_{1} \, - \, c_{1 , N} |$. Now,
\begin{align}
| c_{1} \, - \, c_{1 , N} | &= \ 
\big| \int_{0}^{1} D ( w - w_{N})(s) \, \frac{1}{K(s)} \, ds  \big| 
\ \le  \int_{0}^{1} \big| D ( w - w_{N})(s) \big| \, \frac{1}{K(s)} \, ds    \label{eqDs1} \\
&\le \ 
\int_{0}^{1} \omega^{( (-(\alpha - \beta) + 1)/2 \, , \, (-\beta + 1)/2)}(s) \, \big| D ( w - w_{N})(s) \big| \, 
   \omega^{( (\alpha - \beta - 1)/2 \, , \, (\beta - 1)/2)}(s) \, \frac{1}{K(s)} \, ds  \nonumber \\
&\le \ 
\| D ( w - w_{N} ) \|_{\omega^{( -(\alpha - \beta) + 1 \, , \, -\beta + 1)}} \, 
\bigg\| \frac{1}{K}  \bigg\|_{\omega^{( (\alpha - \beta - 1) \, , \, (\beta - 1))}}    \nonumber   \\
&\le C \, (N + 2)^{- (\alpha-1) } (N(N+\alpha))^{-t \,/ \,2}   \, \| f \|_{t,\,\omega^{(\beta  \, , \, \alpha - \beta )}}  
\ \ \mbox{(using \eqref{error:Dw})} \, .  \label{eqDs2}
\end{align}

In case 
$D(\frac{1}{K}) \in L^{2}_{\omega^{(\alpha - \beta \, , \, \beta)}}(\mrI)$
,
\begin{align}
| c_{1} \, - \, c_{1 , N} | &= \ 
\big| \int_{0}^{1} ( w - w_{N})(s) \, D \big( \frac{1}{K(s)} \big) \, ds  \big|    \label{ester1}  \\
&= \ \int_{0}^{1} \omega^{( -(\alpha - \beta)/2 \, , \, - \beta/2)}(s) \,  ( w - w_{N})(s)  \, 
   \omega^{( (\alpha - \beta)/2 \, , \, \beta/2)}(s) \, D \big( \frac{1}{K(s)} \big) \, ds  \big|  \nonumber \\
&\le \ \| w - w_{N} \|_{  \omega^{( - (\alpha - \beta) \, , \, - \beta)} } \,
 \bigg\| D \big( \frac{1}{K} \big) \bigg\|_{ \omega^{( (\alpha - \beta) \, , \, \beta)} }  \nonumber \\
&\le C \, (N + 2)^{- \alpha } (N(N+\alpha))^{-t \,/ \,2}   \, \| f \|_{t,\,\omega^{(\beta  \, , \, \alpha - \beta )}}
\ \ \mbox{(using \eqref{apertg2})} \, .   \label{ester2} 
\end{align}

We have the following error estimates for $u \, - \, u_{N}$.

\begin{theorem} \label{thcvg4u}
For $f \in H^{t}_{\omega^{(\beta , \alpha - \beta)}}(\mrI)$,  
$t \ge 0$,
then for $\epsilon_{1}, \epsilon_{2} > 0$
there exists $C > 0$ (independent of $N$ and $\alpha$) such that
\begin{align}
\| u \, - \, u_{N} \|_{L^{\infty}} \, + \, \|u-u_N\|_{\omega^{(-1 + \epsilon_{1} \, , \, -1 + \epsilon_{2})}} &\leq \
C \, (N + 2)^{- (\alpha-1) } (N(N+\alpha))^{-t \,/ \,2}   \, \| f \|_{t,\,\omega^{(\beta  \, , \, \alpha - \beta )}}   \, ,  \label{est4u} \\
\  \ \mbox{ and, } \mbox{if 
$D(\frac{1}{K}) \in L^{2}_{\omega^{(\alpha - \beta \, , \, \beta)}}(\mrI)$
},  \nonumber \\
\|u-u_N\|_{\omega^{(-(\alpha-\beta),-\beta)}} &\le  \ C \, 
(N + 2)^{- \alpha } (N(N+\alpha))^{-t \,/ \,2}   \, \| f \|_{t,\,\omega^{(\beta  \, , \, \alpha - \beta )}}  \, .  \label{est4u2} 
\end{align} 
\end{theorem}
\textbf{Proof}: From \eqref{ans4u}, we have using \eqref{eqDs1} and \eqref{eqDs2}
\begin{align}
u(x) \, - \, u_{N}(x)  &= \  (C_{1} \, - \, C_{1 , N}) \int_{0}^{x}  \frac{(1 - s)^{\alpha - \beta -1} \, s^{\beta - 1}}{K(s)} \, ds 
\ + \   \int_{0}^{x}  \frac{D (w \, - \, w_{N})(s)}{K(s)} \, ds  \,    ,   \label{expumuN} \\
\Rightarrow \ \ \| u \, - \, u_{N} \|_{L^{\infty}} &\le \ C \, \big| C_{1} \, - \, C_{1 , N} \big| \ + \  \int_{0}^{1} \big| D (w \, - \, w_{N})(s) \big| \, \frac{1}{K(s)} \, ds  \, , \nonumber \\
&\le \ C 
(N + 2)^{- (\alpha-1) } (N(N+\alpha))^{-t \,/ \,2}   \, \| f \|_{t,\,\omega^{(\beta  \, , \, \alpha - \beta )}}  \, .  \label{neqn2}
\end{align}
Then, from \eqref{neqn2} and 
\begin{align*}
\|u-u_N\|^2_{\omega^{(-1 + \epsilon_{1} \, , \, -1 + \epsilon_{2})}}
&=\int_0^1 \, \omega^{(-1 + \epsilon_{1} \, , \, -1 + \epsilon_{2}}(x) \, (u-u_N)^2(x)  \, dx \\
&\leq   \, \| u-u_N \|^2_{L^\infty} \, \int_0^1 \, \omega^{(-1 + \epsilon_{1} \, , \, -1 + \epsilon_{2})}(x) \, dx
\leq C  \, \| u-u_N \|^2_{L^\infty} \, , 
\end{align*}
we obtain \eqref{est4u}.

For $D(\frac{1}{K}) \in L^{2}_{\omega^{(\alpha - \beta \, , \, \beta)}}(\mrI)$, we apply integration by parts to \eqref{expumuN} to obtain
\begin{align}
u(x)-u_N(x)&=\big(C_1-C_{1,N}\big)\int_0^x\frac{(1-s)^{\alpha-\beta-1}s^{\beta-1}}{K(s)}ds\nonumber\\
& \quad \quad \quad \quad +\frac{w(x)-w_N(x)}{K(x)}+\int_0^x\big(w(s)-w_N(s)\big)\frac{K'(s)}{K^2(s)}ds. \nonumber
\end{align}

Therefore,
\begin{equation}
\|u-u_N\|_{\omega^{(-(\alpha-\beta),-\beta)}}^{2}  \ \le \ I_{1} \ + \ I_{2} \ + I_{3} \, , \quad \mbox{ where }
\label{uu101}
\end{equation}
\begin{align}
I_{1} &= \ 3 \, (C_{1} \, - \, C_{1 , N})^{2} \, \int_{0}^{1}  \omega^{(-(\alpha-\beta),-\beta)}(x) 
\big( \int_0^1 \frac{(1-s)^{\alpha-\beta-1}s^{\beta-1}}{K(s)} ds \big)^{2} \, dx  \,    \nonumber \\
&\le C \, (C_{1} \, - \, C_{1 , N})^{2} \, , 
\label{uu102}
\end{align}
\begin{equation}
I_{2} \ = \ 3 \, \frac{1}{K_{m}^{2}} \, \| w \, - \, w_{N} \|_{\omega^{(-(\alpha-\beta),-\beta)}}^{2} \, ,  \quad \hspace{2.1in} \quad \mbox{ } 
\label{uu103}
\end{equation}
\begin{align}
I_{3} &= \ 3 \, \int_{0}^{1} \omega^{(-(\alpha-\beta) , -\beta)}(x) \big( \int_{0}^{1}
 \omega^{(-\frac{(\alpha-\beta)}{2} , -\frac{\beta}{2})}(s) \, ( w(s) \, - \, w_{N}(s) ) \, 
  \omega^{(\frac{(\alpha-\beta)}{2} , \frac{\beta}{2})}(s) \, \frac{K'(s)}{K^{2}(s)} \, ds \big)^{2} \, dx  \nonumber \\
&\le \ \frac{3}{K_{m}^{2}} \, \int_{0}^{1}    \omega^{(-(\alpha-\beta) , -\beta)}(x) 
\big(  \int_{0}^{1}
 \omega^{(-(\alpha-\beta) , -\beta)}(s) \, ( w(s) \, - \, w_{N}(s) ) \, ds \big) \, 
 \big( \int_{0}^{1}  \omega^{(\alpha-\beta , \beta)}(s) \, \frac{K'(s)}{K^{2}(s)} \, ds \big)  \, dx \nonumber \\
&\le \ C \,   \| w \, - \, w_{N} \|_{\omega^{(-(\alpha-\beta),-\beta)}}^{2} \, .   \label{uu104}
\end{align}
Combining \eqref{uu101}-\eqref{uu104} with \eqref{ester2} and \eqref{apertg2} we obtain \eqref{est4u2}. \\
\mbox{ } \hfill \qed

We conclude this section with an error bound for $D (u - u_{N})$.

\begin{lemma} \label{thcvg4u2}
For $f \in H^{t}_{\omega^{(\beta , \alpha - \beta)}}(\mrI)$,
$t \ge 0$,
then
there exists $C > 0$ (independent of $N$ and $\alpha$) such that
\begin{equation}
\| D ( u - u_{N} ) \|_{\omega^{( -(\alpha - \beta) + 1 \, , \, -\beta + 1)}} \,  \ \le \ C 
(N + 2)^{- (\alpha-1) } (N(N+\alpha))^{-t \,/ \,2}   \, \| f \|_{t,\,\omega^{(\beta  \, , \, \alpha - \beta )}}   \, .  \label{est4u3}
\end{equation}
\end{lemma}
\textbf{Proof}: From \eqref{expumuN} it follows that
\begin{align*}
D(u \, - \, u_{N}) &= \ \frac{1}{K(x)} \left( (C_{1} \, - C_{1,N}) \, (1 - x)^{\alpha - \beta - 1} \, x^{\beta - 1} \ + \ 
 D(w \, - \, w_{N})  \right)  \, . 
\end{align*}
Thus,
\begin{align*}
\| D ( u - u_{N} ) \|_{\omega^{( -(\alpha - \beta) + 1 \, , \, -\beta + 1)}} &\le \ 
   \frac{1}{K_{m}} \bigg(   \, |C_{1} \, - C_{1,N} | \, 
   \| (1 - x)^{(\alpha - \beta - 1)/2} \, x^{(\beta - 1)/2} \|   \\
   &   \quad \quad  \quad \quad  \quad \quad + \ 
  \| D(w \, - \, w_{N})  \|_{\omega^{( -(\alpha - \beta) + 1 \, , \, -\beta + 1)}} \bigg)  \\
&\le \   C 
(N + 2)^{- (\alpha-1) } (N(N+\alpha))^{-t \,/ \,2}   \, \| f \|_{t,\,\omega^{(\beta  \, , \, \alpha - \beta )}} \, ,
\end{align*}
where in the last step we have used \eqref{eqDs2} and \eqref{error:Dw}), and the fact that $(\alpha - \beta - 1)$ and
$(\beta - 1) \, > \, -1$. \\
\mbox{ } \hfill \qed

\setcounter{equation}{0}
\setcounter{figure}{0}
\setcounter{table}{0}
\setcounter{theorem}{0}
\setcounter{lemma}{0}
\setcounter{corollary}{0}

\section{Numerical experiments}
\label{secNum}
In this section we present three numerical experiments to demonstrate our approximation scheme, and to compare
the experimental rate of convergence of the approximation with the theoretically predicated rate. 

\textbf{Numerical example}. Let $K(x)=1/(1+x^\gamma)$ and 
\begin{equation}
f(x) \ = \ -r\frac{x^{1-\alpha}}{\Gamma(2-\alpha)}+(1-r)\frac{(1-x)^{1-\alpha}}{\Gamma(2-\alpha)}.
\label{deffex}
\end{equation}
Then the solution $u(x)$ is given by (\ref{ans4u}) where 
\begin{align*}
w(x) \, = \, x  \, - \, C x^{\beta}{}_2F_1(-(\alpha-\beta-1),\beta;\beta+1,x),~~C = {}_2F_1(-(\alpha-\beta-1),\beta;\beta+1,1)^{-1} \, ,
\end{align*}
and ${}_2F_1(a,b ; c,x)$ denotes the Gaussian three parameter hypergeometric function.

In order to determine the theoretical rate of convergence for $\|u - u_N\|_{\omega^{(-(\alpha-\beta)  , -\beta)}}$,
$\|u - u_N\|_{L^\infty}$ and
$\|D(u - u_N)\|_{\omega^{(-(\alpha-\beta) + 1 \, , \, -\beta + 1)}}$ from \eqref{est4u}, \eqref{est4u2}, and \eqref{est4u3}, respectively,
we need to determine the largest value for $t$ such that $f(x) \in H^{t}_{\omega^{(\beta  \, , \, \alpha - \beta )} }(\mrI)$.
The most singular terms for $f(x)$ in \eqref{deffex} are $x^{1 - \alpha}$ and $(1 - x)^{1 - \alpha}$. 
Using Lemma \ref{lmausp1} (in the Appendix) we have that
$x^{1 - \alpha} \in H^{t}_{\omega^{(\beta  \, , \, \alpha - \beta )} }(\mrI)$, for $t <  \ 3 - \alpha -  \beta$, and
$(1 - x)^{1 - \alpha} \in H^{t}_{\omega^{(\beta  \, , \, \alpha - \beta )} }(\mrI)$, for $t <  \ 3 - \alpha -  (\alpha - \beta)$.

Then, for Experiment 1 ($\alpha = 1.60$, $\beta = 0.80$) 
$f(x) \in H^{t}_{\omega^{(\beta  \, , \, \alpha - \beta )} \,  }(\mrI)$ for $t \, < \, 3 - \alpha - \max\{\alpha - \beta \, , \beta\} = 0.90$, 
which leads to 
theoretical asymptotic convergence rates of 
$\|u - u_N\|_{\omega^{(-(\alpha-\beta)  , -\beta)}} \sim N^{-2.20}$ (using \eqref{est4u2}), 
$\|u - u_N\|_{L^\infty} \sim N^{-1.20}$ (using \eqref{est4u}) and
$\|D(u - u_N)\|_{\omega^{(-(\alpha-\beta) + 1 \, , \, -\beta + 1)}} \sim N^{-1.20}$ (using \eqref{est4u3}).

Assuming that $\| \xi - \xi_{N} \|_{L_{\rho}} \sim  \, N^{- \kappa}$, the experimental convergence rate
is calculated using
\[
  \kappa \approx \frac{\log ( \| \xi - \xi_{N_{1}} \|_{L_{\rho}}  /  \| \xi - \xi_{N_{2}} \|_{L_{\rho}} )}{ \log (N_{2} / N_{1})} \, .
\]

\textbf{Experiment 1}. In this experiment we select $\alpha = 1.60$, $r = 0.50$, $\beta = 0.80$ and $\gamma=0.80$, which leads to
$f(x) \in H^{t}_{\omega^{(\beta  \, , \, \alpha - \beta )} \,  }(\mrI)$ for $t \, < \, 0.60$.

\begin{table}[H]
	\setlength{\abovecaptionskip}{0pt}
	\centering
	\caption{Convergence properties of Experiment 1.}	\label{1,2}
	\vspace{0.5em}	
	\begin{tabular}{ccccccc}
		\hline
		$N$&$\|u-u_N\|_{\omega^{(-(\alpha-\beta),-\beta)}}$ &$\kappa$& $\|D(u-u_N)\|_{\omega^{(-(\alpha-\beta)+1,-\beta+1)}}$&$\kappa$&$\|u-u_N\|_{L^\infty(\mrI)}$&$\kappa$\\
		\cline{1-7}
		16&	4.87E-04&		&1.15E-02	&	&4.41E-04	&\\
		20&	3.09E-04&	2.15&	8.93E-03&	1.18&	2.95E-04&	1.89\\
		24&	2.12E-04&	2.16&	7.26E-03&	1.18&	2.00E-04&	2.23\\
		28&	1.54E-04&	2.16&	6.09E-03&	1.19&	1.54E-04&	1.78\\
		32&	1.16E-04&	2.16&	5.22E-03&	1.19&	1.21E-04&	1.86\\
		36&	9.08E-05&	2.16&	4.56E-03&	1.19&	9.46E-05&	2.14\\
		\hline
		Pred.&    &2.20&     &  1.20&     &  1.20  \\
		\hline	
	\end{tabular}
\end{table}
 
 \begin{figure}[H]
 	\setlength{\abovecaptionskip}{0pt}
 	\centering
 	\includegraphics[width=2.3in,height=2.3in]{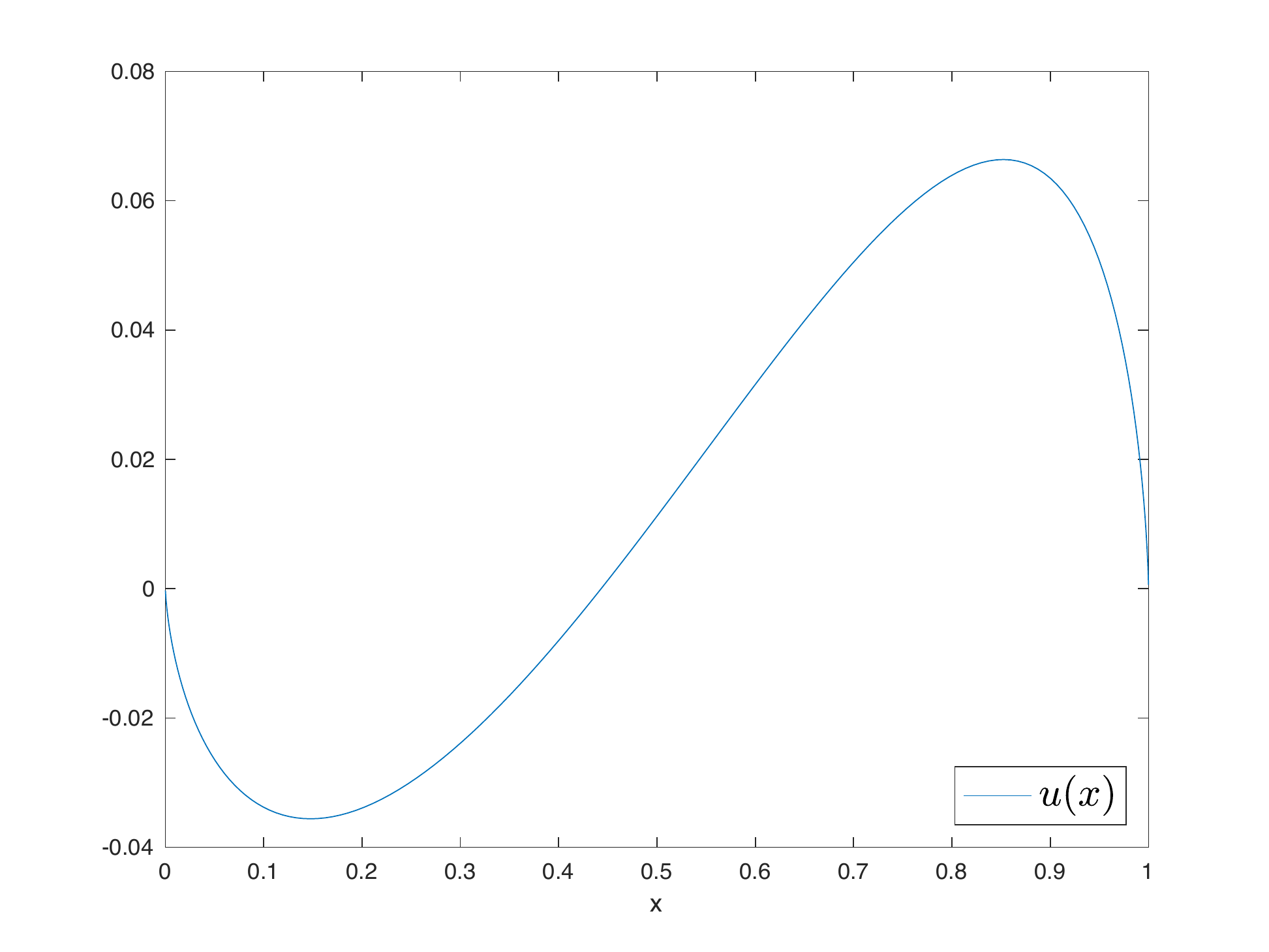}
 	\includegraphics[width=2.3in,height=2.3in]{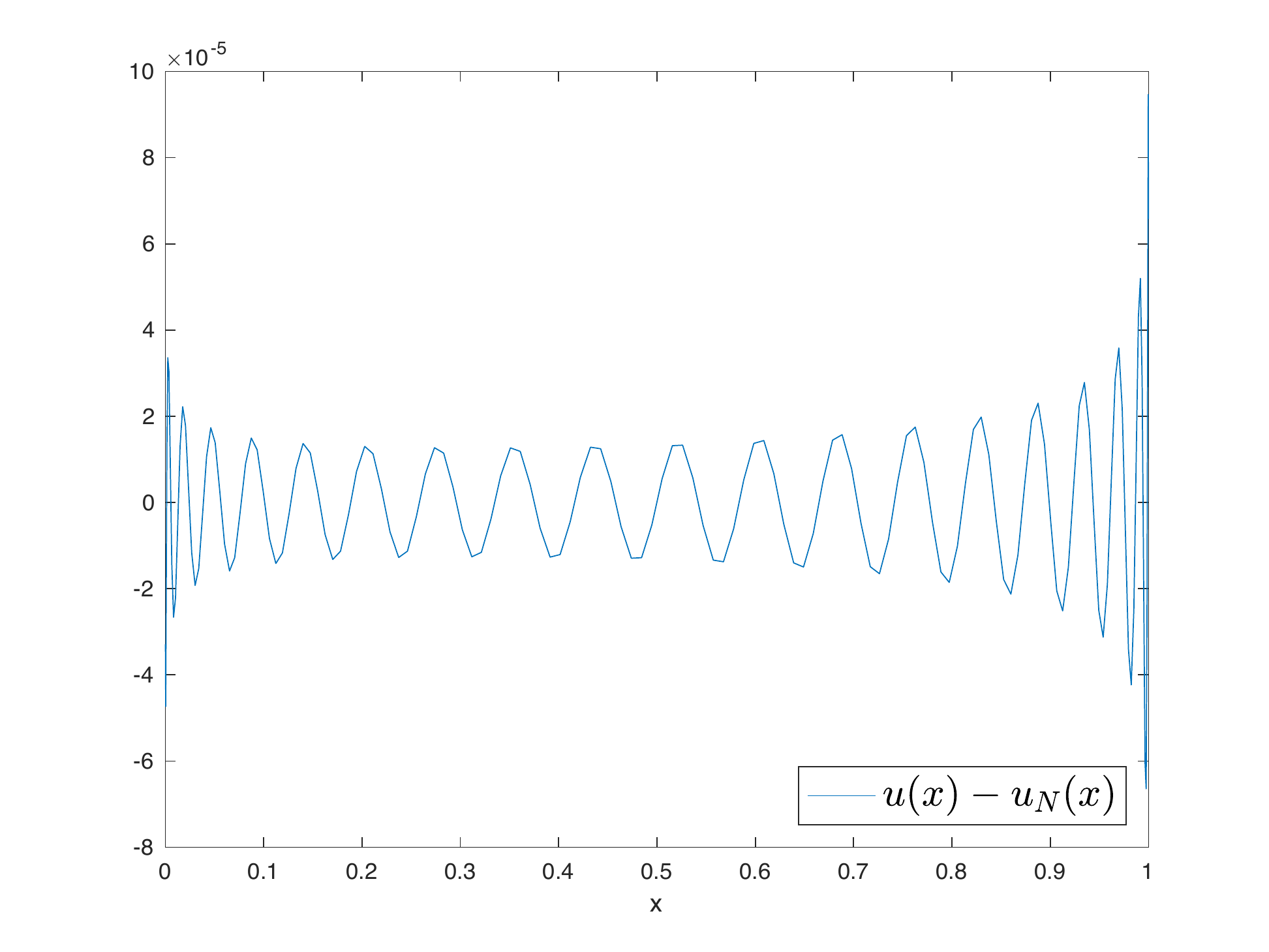}
 	\caption{The plots of (left) $u(x)$ and (right) $u(x)-u_N(x)$}
 \end{figure}

\textbf{Experiment 2}. In this experiment, we take $\alpha = 1.30$, $r = 0.63$, $\beta = 0.50$ and $\gamma=0.80$. The above analysis gives that $f(x) \in H^{t}_{\omega^{(\beta  \, , \, \alpha - \beta )} \,  }(\mrI)$ for $t \, < \,0.90$.
The corresponding theoretical asymptotic convergence rates are 
$\|u - u_N\|_{\omega^{(-(\alpha-\beta)  , -\beta)}} \sim N^{-2.20}$, 
$\|u - u_N\|_{L^\infty} \sim N^{-1.20}$ and
$\|D(u - u_N)\|_{\omega^{(-(\alpha-\beta) + 1 \, , \, -\beta + 1)}} \sim N^{-1.20}$.
\begin{table}[H]
	\setlength{\abovecaptionskip}{0pt}
	\centering
	\caption{Convergence properties of Experiment 2.}	\label{1,1}
	\vspace{0.5em}	
	\begin{tabular}{ccccccc}
		\hline
		$N$&$\|u-u_N\|_{\omega^{(-(\alpha-\beta),-\beta)}}$ &$\kappa$& $\|D(u-u_N)\|_{\omega^{(-(\alpha-\beta)+1,-\beta+1)}}$&$\kappa$&$\|u-u_N\|_{L^\infty(\mrI)}$&$\kappa$\\
		\cline{1-7}
	16&	4.57E-04&	&	1.07E-02&		&4.41E-04&	\\
	20&	2.88E-04&	2.20&	8.29E-03&	1.21&	2.95E-04&	1.90\\
	24&	1.96E-04&	2.19&	6.71E-03&	1.21&	2.00E-04&	2.24\\
	28&	1.42E-04&	2.19&	5.61E-03&	1.21&	1.53E-04&	1.78\\
	32&	1.07E-04&	2.19&	4.80E-03&	1.21&	1.20E-04&	1.87\\
	36&	8.33E-05&	2.19&	4.18E-03&	1.21&	9.42E-05&	2.15\\
		\hline
		Pred.&    &2.20&     &  1.20&     &  1.20  \\
		\hline	
	\end{tabular}
\end{table}

\textbf{Experiment 3}. In this experiment we select $\alpha = 1.30$, $r = 0.63$, $\beta = 0.50$ and $\gamma=0.10$, which leads to
$f(x) \in H^{t}_{\omega^{(\beta  \, , \, \alpha - \beta )} \,  }(\mrI)$ for $t \, < \, 0.90$. However, in this case, 
$D(\frac{1}{K}) \not \in L^{2}_{\omega^{(\alpha - \beta \, , \, \beta)}}(\mrI)$ due to the relatively strong singularity of $K(x)$ at $x=0$, which means that (\ref{ester2}) is not applicable. 
Hence we can only apply the bound (\ref{eqDs2}) of $c_1-c_{1,N}$, which leads to the estimate (\ref{est4u}) of $u-u_N$ instead of (\ref{est4u2}), and consequently an estimate for the
convergence rate of $\|u-u_N\|_{\omega^{(-(\alpha-\beta),-\beta)}}$ of $1.20$ by (\ref{est4u}), instead of $2.20$ 
if using (\ref{est4u2}).

\begin{table}[H]
	\setlength{\abovecaptionskip}{0pt}
	\centering
	\caption{Convergence properties of Experiment 3.}	\label{1,3}
	\vspace{0.5em}	
	\begin{tabular}{ccccccc}
		\hline
		$N$&$\|u-u_N\|_{\omega^{(-(\alpha-\beta),-\beta)}}$ &$\kappa$& $\|D(u-u_N)\|_{\omega^{(-(\alpha-\beta)+1,-\beta+1)}}$&$\kappa$&$\|u-u_N\|_{L^\infty(\mrI)}$&$\kappa$\\
		\cline{1-7}
		16&	4.99E-04&		&1.14E-02	&	&5.73E-04	&\\
		20&	3.11E-04&	2.24&	8.74E-03&	1.25&	3.71E-04&	2.05\\
		24&	2.10E-04&	2.24&	7.03E-03&	1.24&	2.63E-04&	1.98\\
		28&	1.51E-04&	2.23&	5.85E-03&	1.24&	1.88E-04&	2.25\\
		32&	1.13E-04&	2.22&	4.99E-03&	1.24&	1.35E-04&	2.59\\
		36&	8.80E-05&	2.22&	4.33E-03&	1.23&	1.07E-04&	2.01\\
		\hline
		Pred.&    &1.20&     &  1.20&     &  1.20  \\
		\hline	
	\end{tabular}
\end{table}


The experimental convergence rates for $\|u-u_N\|_{L^2_{\omega^{(-(\alpha-\beta),-\beta)}}}$ and $\|D(u-u_N)\|_{L^2_{\omega^{(-(\alpha-\beta)+1,-\beta+1)}}}$ are in strong agreement with the theoretically predicted rates for the first two experiments. For Experiment 3, we note that the numerical convergence rate of $\|u-u_N\|_{L^2_{\omega^{(-(\alpha-\beta),-\beta)}}}$ is 2.20, which corresponds to the case for $D(\frac{1}{K}) \in L^{2}_{\omega^{(\alpha - \beta \, , \, \beta)}}(\mrI)$
(even though this is not the case in this experiment). 
Additionally, we remark that the error in the $L^\infty$ norm is difficult to measure accurately due to the oscillatory nature of the error function caused by the polynomial approximation, as illustrated by Figure 1.

\appendix

\setcounter{equation}{0}
\setcounter{figure}{0}
\setcounter{table}{0}
\setcounter{theorem}{0}
\setcounter{lemma}{0}
\setcounter{corollary}{0}

\section{Appendix}
\label{apdxAPP}
In this section we investigate which $H^{s}_{\omega^{(\alpha , \beta)}}(\mrI)$ space $u(x) = x^{\mu}$ lies in.
For brevity of notation, in this section we use 
$H^{s}_{(\alpha , \beta)}(\mrI) \equiv H^{s}_{\omega^{(\alpha , \beta)}}(\mrI)$.

\begin{lemma}  \label{lmausp1}
Let $u(x) \, = \, x^{\mu}$. Then, $u \in H^{s}_{(\alpha , \beta)}(\mrI)$ for $s \, \ge \, 0$ 
satisfying $s  \, < \, 2 \mu + \beta + 1$. \\
\end{lemma}
\textbf{Proof}: Let $\chi(x) \in C^{\infty}[0 , \infty)$ denote the cutoff function satisfying
\[
\chi(x) \ = \ \left\{ \begin{array}{rl}
1 & \mbox{ for } \ 0 < x \le 1/4  \\  0 & \mbox{ for } \   x \ge  3/4 \end{array}  \right. \, ,
\]
and let $\chi_{\delta}(x) \, := \, \chi(\frac{x}{\delta})$, for $\delta > 0$. Note that
\[
\chi_{\delta}(x) \ = \ \left\{ \begin{array}{rl}
1 & \mbox{ for } \ 0 < x \le \delta/4  \\  0 & \mbox{ for } \   x \ge  3 \delta /4 \end{array}  \right. \, ,
\quad \mbox{ and } \quad
\frac{d^{m}}{dx^{m}} \chi_{\delta}(x) \ = \ \left\{ \begin{array}{rl}
0 & \mbox{ for } \ 0 < x <  \delta/4  \\  0 & \mbox{ for } \   x  >  3 \delta /4 \end{array}  \right. \, 
\mbox{ for } m \in \mathbb{N} \, .
\]

For $\delta$ to be determined, let $u \, = \, v + w$ where $v(x) \ = \ \chi_{\delta}(x) \, u(x)$ and
$w(x) \ = \ \big(1 - \chi_{\delta}(x) \big) \, u(x)$. We have that
\[
\left| \frac{d^{m} v(x)}{dx^{m}} \right| \ \le \ C \, \sum_{j = 0}^{m} \delta^{-(m - j)} \, x^{\mu - j} \, , 
\ \mbox{ and  is zero for }  \ x \ > \ 3 \delta / 4 \, .
\]
Thus,
\begin{align*}
\int_{I} (1 - x)^{\alpha + m} \, x^{\beta + m} \, \left| \frac{d^{m} v(x)}{dx^{m}} \right|^{2} \, dx
&\le \  C \, \int_{0}^{3 \delta / 4} \, \sum_{j = 0}^{m} \delta^{-2 (m - j)} \, x^{2 \mu \, - \, 2 j \, + \beta + m} \, dx  \\
&\le \ C \, \sum_{j = 0}^{m} \delta^{-2 (m - j)} \, \delta^{2 \mu \, - \, 2 j \, + \beta + m + 1}
 \ \ \mbox{ provided } 2 \mu \, - \, 2 j \, + \beta + m > -1 \, ,   \\
 &\le \ C \, \delta^{2 \mu \,  + \, \beta + 1 - m} \, , 
\end{align*}
which implies that, for $m \, < \, 2 \mu \,  + \, \beta + 1$, $v \in H^{m}_{(\alpha ,  \beta)}(\mrI)$ and
\be
\| v \|_{H^{m}_{(\alpha ,  \beta)}}^{2} \ \le \ C \, \delta^{2 \mu \,  + \, \beta + 1 - m} \, .
\label{ewst21}
\ee

Next, consider $w(x)$.
\be
\left| \frac{d^{m} w(x)}{dx^{m}} \right| \ \le \ C \, \left( \big(1 - \chi_{\delta}(x) \big) x^{\mu - m} \ + \ 
 \sum_{j = 0}^{m - 1} x^{\mu - j} \,  \frac{d^{m - j}}{dx^{m - j}} \big(1 - \chi_{\delta}(x) \big) \right)  \,  . 
\label{wterm1}
\ee
The first term on the RHS of \eqref{wterm1} vanishes for $x < \delta / 4$, and the second term vanishes for
$x \, < \, \delta/4$ and $x \, > \, 3 \delta/4$. Using this,
\begin{align}
\int_{I} (1 - x)^{\alpha + m} \, x^{\beta + m} \, \left| \frac{d^{m} w(x)}{dx^{m}}  \right|^{2} \, dx
&\le \ C \, \left( \int_{\delta/4}^{1} x^{\beta + m} \, x^{2 \mu \, - \, 2 m} \, dx \ + \ 
\int_{\delta / 4}^{3 \delta / 4} \,  x^{\beta + m} \,  
 \sum_{j = 0}^{m - 1} x^{2 \mu \, - \, 2 j} \,  \delta^{-2 ( m - j)} \,  dx \right)   \nonumber \\
&\le \  C \, \left( \int_{\delta/4}^{1}  x^{2 \mu \, + \beta \, - \, m} \, dx \ + \  
 \sum_{j = 0}^{m - 1} \, \delta^{-2 m \, + \, 2 j} \, \int_{\delta / 4}^{3 \delta / 4} \, x^{2 \mu \, + \,  \beta \, + \, m \, - \, 2 j} \,  dx \right)
 \nonumber  \\
&
\le \ \left\{ \begin{array}{rl}
           C \, \left( 1 \, + \, \delta^{2 \mu \, + \, \beta \, - \, m \, + \, 1}  \right) & \mbox{ if }  2 \mu \, + \, \beta \, - \, m \, \neq \, -1 \\
           C \, \left( 1 \, + \, | \log \delta |  \right) & \mbox{ if }  2 \mu \, + \, \beta \, - \, m \, = \, -1 
           \end{array} \right. \, .  \label{ewst23}
\end{align}
Hence, for $n \, > \, 2 \mu \, + \, \beta \, + \, 1$
\be
   \| w \|_{H^{n}_{(\alpha ,  \beta)}}^{2} \ \le \ C \, \delta^{2 \mu \, + \, \beta \, + \, 1 \, - \, n} \, .
   \label{ewst24}
\ee
(\textbf{Remark}: For $n \, > \, 2 \mu \, + \, \beta \, + \, 1$ the exponent of $\delta$ in \eqref{ewst24} is negative,
so the $'1'$ term in \eqref{ewst23} is bounded by the $\delta$ term.)

For $0 < t < 1$ we have from \eqref{ewst21} and \eqref{ewst24}
\begin{align}
K(t , u) &= \ \inf_{u \, = \, u_{1} + u_{2}} \left( \| u_{1} \|_{H^{m}_{(\alpha , \beta)}} \ + \ t \, \| u_{1} \|_{H^{n}_{(\alpha , \beta)}} 
 \right)  \label{ewst251}  \\
 &\le \ \| v \|_{H^{m}_{(\alpha , \beta)}} \ + \ t \, \| w \|_{H^{n}_{(\alpha , \beta)}}   \nonumber \\
 &\le \ C \, \left( \delta^{(2 \mu \,  + \, \beta + 1 - m) / 2} \ + \ t \, \delta^{(2 \mu \, + \, \beta \, + \, 1 \, - \, n) / 2} \right) \, .
 \label{ewst252}
\end{align} 

Setting $\delta \ = \ t^{2 / (n - m)}$ leads to  $K(t , u) \, \le \, C \, t^{(2 \mu \,  + \, \beta + 1 - m) / (n - m)}$.

Recall that
\be
\| u \|^{2}_{[H^{m}_{(\alpha ,  \beta)} \, , \, H^{n}_{(\alpha ,  \beta)}]_{\theta , 2}} \ = \ 
\int_{0}^{\infty} \, t^{-2 \theta} \left( K(t , u) \right)^{2} \, \frac{dt}{t} \, .
\label{ewst26}
\ee
The larger the value of $\theta$ $(0 < \theta < 1)$ in \eqref{ewst26} such that the integral is finite, 
the ``nicer'' (i.e., more regular) is the function $u$. Hence from 
\eqref{ewst26}, we are interested in the integrand about $t = 0$. We have trivially that for $u_{1} = u$, $u_{2} = 0$
in \eqref{ewst251} that $K(t , u) \, \le \, \| u \|_{H^{m}_{(\alpha ,  \beta)}}  \, \le \, C$. Hence it follows that
\[
  K(t , u) \, \le \, \left\{ \begin{array}{ll}
  C \, t^{(2 \mu \,  + \, \beta + 1 - m) / (n - m)} & \mbox{ for } 0 < t < 1 \\
  C &  \mbox{ for } t \ge 1 
  \end{array} \right.   \, .
\]  

Using \eqref{ewst26},
\begin{align}
\| u \|^{2}_{[H^{m}_{(\alpha ,  \beta)} \, , \, H^{n}_{(\alpha ,  \beta)}]_{\theta , 2}} &\le \
\int_{0}^{1} C \, t^{-2 \theta \, - \, 1 \, + \, 2 (2 \mu \,  + \, \beta + 1 - m) / (n - m)} \, dt \ + \ 
\int_{1}^{\infty}  C \, t^{-2 \theta \, - \, 1} \, dt \ < \ \infty \, , \nonumber  \\
\mbox{ if } \ \ \theta &< \ (2 \mu \,  + \, \beta + 1 - m) / (n - m) \, .  \label{ewst27}
\end{align}
For $s \ = \ (1 - \theta) m \ + \ \theta n \ = \ m \ + \ \theta (n - m)$, then  
$s \ < \ 2 \mu \,  + \, \beta + 1$ using  \eqref{ewst27} \, .

Hence we can conclude that $u(x) \ = \ x^{\mu} \, \in \, H^{s}_{(\alpha ,  \beta)}(\mrI)$ for $s < 2 \mu \,  + \, \beta + 1$. \\
\mbox{  } \hfill \qed

\end{document}